
\input amstex.tex
\documentstyle{amsppt}
\hoffset=-1in \advance \hoffset by 22mm \pagewidth{138mm}
\voffset=-1in \advance \voffset by 36mm \pageheight{185mm}
\NoRunningHeads
\topmatter
\title
{Algebraic surfaces with quotient singularities -
including some discussion on automorphisms and fundamental groups}
\endtitle
\author
J. Keum and D. -Q. Zhang
\endauthor
\endtopmatter

\document

\head
Introduction
\endhead

We work over the complex numbers field ${\bold C}$.
In the present survey, we report some recent progress on the study of varieties
with mild singularities like log terminal singularities
(which are just quotient singularities in the case of dimension 2; see [KMM]).
Singularities appear naturally in many ways. The minimal model program
developed by Mori et al shows that a minimal model will inevitably have some terminal
singularities [KMo]. Also the degenerate fibres of a family of varieties will have some
singularities.

\par
We first follow Iitaka's strategy to divide (singular) varieties $Y$ according to
the logarithmic Kodaira dimension $\kappa(Y^0)$ of the smooth locus $Y^0$ of $Y$.
One key result {\bf (2.3)} says that for a relatively minimal log terminal surface $Y$
we have either nef $K_Y$ or dominance of $Y^0$ by an affine-ruled surface.
It is conjectured to be true for any dimension [KMc].

\par
In smooth projective surfaces of general type case, we have Miyaoka-Yau inequality
$c_1^2 \le 3c_2$
and Noether inequalities: $p_g \le (1/2)c_1^2 + 2, \,\,$
$c_1^2 \ge (1/5)c_2 - (36/5)$. Similar inequalities are given for $Y^0$ in Section 4;
these will give effective restriction on the region where non-complete
algebraic surfaces of general type exist.

\par
In Kodaira dimension zero case, an interesting conjecture {\bf (3.12)}
(which is certainly true when $Y$ is smooth projective by the classification theory)
claims that for a relatively minimal and log terminal surface $Y$
of Kodaira dimension $\kappa(Y^0) = 0$, one has either
$\pi_1(Y^0)$ finite, or an etale cover $Z^0 \rightarrow Y^0$
where $Z^0$ is the complement of a finite set in an abelian surface $Z$.
Some partial answers to {\bf (3.12)} are given in Section 3.

\par
The topology of $Y^0$ is also very interesting. We still do not know whether
$\pi_1$ of the complement of a plane curve is always residually finite or not.
Conjecture {\bf (2.4)} proposed in [Z7] claims that the smooth locus of a log terminal
Fano variety has finite topological fundamental group. This is confirmed when
the dimension is two and now there are three proofs: [GZ1, 2] (using Lefschetz hyperplane
section theorem and Van Kampen theorem), [KMc] (via rational connectivity), [FKL] (geometric).

\par \vskip 1pc
The other interesting topic covered is the automorphism groups.
Recent progress in K3 surface case is treated in Section 5. For
generic rational surface $X$ of degree $\le 5$, it is classically
known that $|\text{\rm Aut}(X)|$ divides $5!$. However, when $Y$ is a log
terminal del Pezzo singular surface of Picard number 1, it is
very often that $\text{\rm Aut}(Y)$ contains ${\bold Z}/(p)$ for all prime
$p \ge 5$ (see [Z9] or {\bf (6.2)}).

\head
Terminology and Notation
\endhead

\par \noindent
{\bf (1).} For a variety $V$ we denote by $V^0 = V -$\text{\rm Sing} \, $V$ the smooth
locus.

\par \noindent
{\bf (2).} A $(-n)$-curve $C$ on a smooth surface is a smooth rational curve
with $C^2 = -n$.

\par \noindent
{\bf (3).} For a divisor $D$, we denote by $\#D$ the number of irreducible
components of Supp $D$.

\par \noindent
{\bf (4).} For a variety $V$, the $e(V)$ is the Euler number.

\head
Section 1. Preliminaries
\endhead

\par \noindent
{\bf (1.1).} Let $V^0$ be a nonsingular variety and let
$V$ be a {\it smooth completion} of $V^0$, i.e., $V$ is nonsingular projective and
$D := V \setminus V^0$ is a divisor with simple normal crossings.
If $H^0(V, m(K_V+D)) = 0$ for all $m \ge 1$, we define the
{\it Kodaira (logarithmic) dimension} $\kappa(V^0) = -\infty$.
Otherwise, $|m(K_V+D)|$ gives rise to a rational map $\varphi_m$
for some $m$ and we define the {\it Kodaira dimension}
$\kappa(V^0)$ as the maximum of $\dim(\varphi_m(V^0))$.

\par
The Kodaira dimension of $V^0$ does not depend on the choice of
the completion $V$ [I3, \S 11.2]. Also
$\kappa(V^0)$ takes value in $\{-\infty, 0, 1, \dots, \dim V^0\}$.

\par
$V^0$ is of {\it general type} if $\kappa(V^0) = \dim V^0$.
\par \vskip 1pc
$p_g(V^0) = h^0(V, K_V+D)$ is called the {\it logarithmic geometric genus}
which does not depend on the choice of the completion $V$ [ibid.].

\par \vskip 1pc \noindent
{\bf (1.2).}

\par \noindent
{\bf (a)} Let $G \subseteq GL_2({\bold C})$ be a non-trivial finite group
with no reflection elements. Then ${\bold C}^2/G$ has a unique
singularity at $\overline{O}$
(the image of the origin of the affine plane ${\bold C}^2$).
A singularity $Q$ of a normal surface $Y$ is a {\it quotient singularity}
if locally the germ $(V, Q)$ is analytically isomorphic to
$({\bold C}^2/G, \overline{O})$
for some $G$. Quotient singularities are classified in [Br, Satz 2.11].

\par \noindent
{\bf (b)} When $\dim Y = 2$, the $Q$ in $Y$ is a quotient singularity if
and only if it is a log terminal singularity [Ka2, Cor 1.9].

\par
$Q$ is a {\it Du Val} (or {\it rational double}, or {\it Dynkin type} $ADE$,
or {\it canonical}, or {\it rational Gorenstein}
in other notation) singularity if $G \subseteq SL_2({\bold C})$
(see [Du], [Re1]).

\par \vskip 1pc
In {\bf (1.3) - (1.6)} below, we assume that $Y$ is
a normal projective surface with at worst quotient singularities.

\par \vskip 1pc \noindent
{\bf (1.3).} Let $f: {\widetilde Y} \rightarrow Y$ be
the minimal resolution and $D$ the exceptional divisor.
We can write $f^*K_Y = K_{\widetilde Y} + D^*$ where $D^*$ is an effective
{\bf Q}-divisor with support in $D$. Write $D = \sum_{i=1}^n D_i$
with irreducible $D_i$ and $D^* = \sum_{i=1}^n d_i D_i$.

\par \vskip 1pc \noindent
{\bf Lemma.}

\par \noindent
(1) Each $D_i$ is a $(-n_i)$-curve for some $n_i \ge 2$.

\par \noindent
(2) $0 \le d_i < 1$.

\par \noindent
(3) $d_i = 0$ holds if and only if the connected component of $D$
containing $D_i$ is contracted to a Du Val singularity on $Y$
(i.e., $f(D_i)$ is a Du Val singularity).

\par \vskip 1pc
{\it Proof.} (2) follows from the fact that a quotient singularity
is just a log terminal singularity [Ka2, Cor 1.9]. For (1) and (3), see
[Br, Satz 2.11] and [Ar1, Theorem 2.7].

\par \vskip 1pc \noindent
{\bf (1.4).} $Y$ is a {\it log del Pezzo} surface
if the anti-canonical divisor $-K_Y$ is ${\bold Q}$-ample.
$Y$ is a {\it Gorenstein (log) del Pezzo}
surface if further $Y$ has at worst Du Val singularities
(see [Mi, Ch II, 5.1], [MZ1]).

\par
A normal variety $V$ is {\it Fano} if $-K_V$ is {\bf Q}-ample.

\par
A log del Pezzo surface is nothing but a log terminal
Fano surface, and a Gorenstein del Pezzo surface is nothing
but a canonical Fano surface.

\par \vskip 1pc \noindent
{\bf (1.5).} $Y$ is a {\it log Enriques surface} if the
irregularity $q(Y) = h^1(Y, {\Cal O}_Y) = 0$ and if $mK_Y \sim 0$
(linear equivalence) for some positive integer $m$. The smallest
$m$ is called the index of $Y$ and denoted by $I(Y)$ [Z4, Part I,
Definition 1.1].

\par
A log Enriques surface of index 1 is nothing but a K3 surface
possibly with Du Val singularities. A non-rational log Enriques surface
is of index 2 if and only if it is an Enriques surface possibly with
Du Val singularities. The case of $Y$ with a unique singularity
is classified by Tsunoda [Ts, Proposition 2.2] (see also [Z4, Part I, Proposition 1.6]).

\par \vskip 1pc \noindent
{\bf Proposition.} {\it Let $Y$ be a rational log Enriques surface with $\#(\text{\rm Sing} \, Y) = 1$.
Then $I(Y) = 2$ and the unique singularity is of type $(1/4n)(1, 2n-1)$
for some $n \ge 1$.}

\par \vskip 1pc \noindent
{\bf (1.6).} The surface $Y$ is {\it relatively minimal}
if for every curve $C$, we have either $K_Y . C \ge 0$ or $C^2 \ge 0$.

\par
Suppose that $Y$ is not relatively minimal.
Then there is a curve $C$ such that $K_Y . C < 0$
and $C^2 < 0$. By [MT1, Lemma 1.7 (2)], we see that there is a contraction
$Y \rightarrow Z$ of the curve $C$ to a smooth or quotient singularity
such that the Picard number $\rho(Z) = \rho(Y) - 1$.
So every projective surface with at worst quotient singularities has a
relatively minimal model.

\par \vskip 1pc
$Y$ is {\it strongly minimal} if it is relatively minimal and
if there is no curve $C$ with $C^2 < 0$ and $C . K_Y = 0$ [Mi, Ch II, {\bf (4.9)}].

\par \vskip 1pc \noindent
{\bf (1.7).} A smooth projective rational surface $X$
is a {\it Coble} surface if $|-K_X| = \emptyset$
while $|-2K_X| \ne \emptyset$. A Coble surface is
{\it terminal} if it is not the image of any birational but
not biregular morphism of Coble surface.

\par
Coble surfaces are classified in [DZ]. Here is an example.
Let $Z$ be a rational elliptic surface with a multiplicity-2 fibre $F_0$
and a non-multiple fibre $F_1$ of type $I_n$ (see [CD] for classification of $Z$).
Let $X \rightarrow Z$ be the blow up of all $n$ intersection points in $F_1$.
Then $X$ is a terminal Coble surface. Coble surfaces and log Enriques
surfaces are closely related.

\par \vskip 1pc \noindent
{\bf Proposition} [DZ, Proposition 6.4].

\par \noindent
(1) {\it The minimal resolution $X$ of a rational log Enriques
surface $Y$ of index $2$ is a Coble surface with $h^0(X, -2K_X) = 1$
and the only member $D$ in $|-2K_X|$ is reduced and a disjoint union of
$D_i$, where $D_i$ is either a single $(-4)$-curve or a linear
chain with the dual graph below (each $D_i$ is contractible
to a singularity of type $(1/4n_i)(1, 2n_i-1)$ with $n_i = \#D_i$):}
$$(-3)--(-2)--\cdots--(-2)--(-3).$$

\par \noindent
{\it If we let $X_{te} \rightarrow X$ be the blow up of all intersection
points in $D$, then $X_{te}$ is a terminal Coble surface
with $h^0(X_{te}, -2K_{X_{te}}) = 1$
and the only member in $|-2K_{X_{te}}|$ is a disjoint union
of $n$ of $(-4)$-curves with $n = \sum_i n_i$.}
\par \noindent
(2) {\it Conversely, a terminal Coble surface $X$ has a unique member
$D$ in $|-2K_X|$, and $D$ is reduced and a disjoint union of $(-4)$-curves.}

\par \vskip 1pc \noindent
{\bf (1.8).} A smooth affine surface $S$ is a {\bf Q}-{\it homology plane}
if $H_i(S, {\bold Q}) = 0$ for all $i > 0$.
Similarly we can define a {\bf Z}-{\it homology plane}
and {\bf Q}-{\it homology plane with quotient singularities}.
The following very important theorem is proved by
Gurjar, Pradeep and Shastri in their papers
[GS], [PS], [GPS] and GPr].

\par \vskip 1pc \noindent
{\bf Theorem 1.9.} {\it A {\bf Q}-homology plane with at worst quotient
singularities is a rational surface.}

\par \vskip 1pc \noindent
{\bf Theorem 1.10} [Mi, Theorem 4.10]. {\it Let $Y$ be a ${\bold Q}$-homology
plane. Let $\nu$ be the number of topologically contractible curves in $Y$.
Then we have:}

\par \noindent
(1) {\it Every topologically contractible curve is isomorphic to the affine line.}

\par \noindent
(2) {\it The Kodaira dimension $\kappa(Y^0) = 2, 1$ or $0, -\infty$
if and only if $\nu = 0,$ finite, $\infty$, respectively.}
(Gurjar and Parameswaran [GPa] have determined the number $\nu$ when
$\kappa(Y^0) = 0$).

\par \noindent
(3) {\it Suppose that $Y$ is a homology plane. Then
$\kappa(Y^0) = 2, 1, -\infty$
if and only if $\nu = 0, 1, \infty$, respectively}
(see {\bf (3.17)}).

\head
Section 2. Normal Algebraic Surfaces $Y$ with Kodaira dimension
$\kappa(Y^0) =-\infty$
and Fano varieties
\endhead

In this section we consider projective varieties with
at worst log terminal singularities and
Kodaira dimension $\kappa(Y^0) = -\infty$,
where $Y^0 = Y - \text{\rm Sing} \, Y$.

\par
The following result is a special case of [MT1, Theorem 2.11].
We will sketch a different and direct proof here by making use of [KMM].

\par \vskip 1pc \noindent
{\bf Theorem 2.1.} {\it Let $Y$ be a relatively minimal surface
with at worst quotient singularities. Then one of the following occurs.}

\par \noindent
(1) {\it The Kodaira dimension $\kappa(Y^0) \ge 0$
and $K_Y$ is numerically effective.}

\par \noindent
(2) {\it $\kappa(Y^0) = -\infty$ and $K_Y$ is not numerically effective.
To be precise, either}

\par \noindent
(2a) {\it $Y^0$ is ruled, i.e. $Y^0$ has a Zariski open set
of the form ${\bold P}^1 \times C$ with a curve $C$, or}

\par \noindent
(2b) {\it $Y$ is a log del Pezzo surface of Picard number $1$.}

\par \vskip 1pc
{\it Proof.} We may assume that $K_Y$ is not nef. Then by
[KMM, Theorems 4-2-1 and 3-2-1], there is an extremal ray
${\bold R}_{> 0} [C]$ with $C$ a rational curve,
and a corresponding morphism $\Phi : Y \rightarrow Z$
with connected fibres such that a curve $E$ is
mapped to a point by $\Phi$ if and only
if the class of $E$ is in ${\bold R}_{> 0}[C]$.
\par \vskip 1pc \noindent
{\bf Case $\dim Z = 2$.} Then $Z$ has at worst log terminal singularities
(= quotient singularities) by [KMM, Proposition 5-1-6]. This contradicts the relative
minimality of $Y$.

\par \vskip 1pc \noindent
{\bf Case $\dim Z = 0$.} Then $\text{\rm Pic}  \,  Y$ is generated over ${\bold Q}$
by $C$ and hence Picard number $\rho(Y) = 1$ and $C$ is ${\bold
Q}$-ample. Since $K_Y . C < 0$, we have $K_Y = aC$ (numerically)
with $a < 0$. So the case (2b) occurs.

\par \vskip 1pc \noindent
{\bf Case $\dim Z = 1$.} Then a general fibre $F$ of $\Phi$ is
${\bold P}^1$ because $F^2 = 0$ and $K_Y . F < 0$ (pull back to
${\widetilde Y}$ and use genus formula). Clearly the case (2a)
occurs. This proves the theorem.

\par \vskip 1pc \noindent
{\bf Theorem 2.2} ([KMc, Cor. 1.6], [Mi, Ch II, Theorems 2.1 and 2.17]).

\par \noindent
{\it Let $Y$ be a log del Pezzo surface.
Then there is a dominant morphism $X^0 \rightarrow Y^0$ such that
$X^0$ is an affine-ruled surface (i.e., $X^0$ contains a Zariski open set
of the form ${\bold A}^1 \times C$ for some curve $C$).}

\par \vskip 1pc
When $Y$ is Gorenstein, Theorem 2.2 was proved in [Z2, Theorem
3.6]; the general case of Theorem 2.2 was proved in a lengthy
book [KMc, Cor 1.6]; Theorems 2.1 and 2.2 together give the proof
of the following result, which is the quotient surface case of
Miyanishi Conjecture.

\par \vskip 1pc \noindent
{\bf Theorem 2.3.} {\it Let $Y$ be a projective surface with at worst quotient
singularities. Suppose that $Y$ is relatively minimal.
Then the following are equivalent.}

\par \noindent
(1) {\it $K_Y$ is not nef.}

\par \noindent
(2)  {\it $\kappa(Y^0) = -\infty$.}

\par \noindent
(3) {\it There is a dominant morphism $X^0 \rightarrow Y^0$
such that $X^0$ is an affine-ruled surface.}

\par \vskip 1pc \noindent
{\it Proof.} The equivalence of (1) and (2) is proved in [MT1, Theorem 2.11].
(2) implies (3) by Theorems 2.1 and 2.2.
Assume (3). Now $\kappa(X^0) = -\infty$ is clear by considering a ruled surface
as a completion of $X^0$ with the boundary equal to the union of
a section (or empty set) and a few fibre components.
Since $\kappa(X^0) \ge \kappa(Y^0)$, (2) follows.

\par \vskip 1pc
Now we turn to the topology of smooth locus of a variety.
We proposed the following in [Z7].

\par \vskip 1pc \noindent
{\bf Conjecture 2.4.} {\it Let $V$ be a Fano variety with
at worst log terminal singularities.
Then the topological fundamental group $\pi_1(V^0)$ is finite.}

\par \vskip 1pc
The affirmative answer to {\bf (2.4)} would imply the following
which was conjectured
in [KZ] and is now a theorem of S. Takayama [Ta].
Indeed, {\bf (2.4)} would imply that $\pi_1(V)$ is finite
and we let $U \rightarrow V$ be the universal cover.
Then $\chi({\Cal O}_U) = n \chi({\Cal O}_V)$, where
$n = |\pi_1(V)|$. The Kawamata-Viehweg vanishing
implies that $\chi({\Cal O}_X) = h^0(X, {\Cal O}_X)$ ($= 1$)
for both $X = U$ and $V$. Hence $n = 1$.

\par \vskip 1pc \noindent
{\bf Theorem 2.5.} {\it Let $V$ be a Fano variety
with at worst log terminal singularities. Then $\pi_1(V) = (1)$.}

\par \vskip 1pc
The result {\bf (2.5)} would also follow from the following conjecture which is still open for dimension 4 or higher. It is proved in 3-fold case by [Ca] and [KoMiMo].
Recently, Graber, Harris and Starr [GHS] have proved that
any complex algebraic variety having a fibration with rationally connected
general fibres and image (or base), is again rationally connected.

\par \vskip 1pc \noindent
{\bf Conjecture 2.6.} {\it Let $V$ be a Fano variety with log terminal singularities.
Then $V$ is rationally connected, i.e., any two general points are connected by
an irreducible rational curve.}

\par \vskip 1pc
Partial answers to {\bf (2.4)} are given in {\bf (2.7)} $\sim$ {\bf (2.10)}.

\par \vskip 1pc \noindent
{\bf (2.7).} When dimension is 2, Conjecture 2.4 was proved in
affirmative by [GZ1, 2]; for a differential geometric proof, see
[FKL]. In [KMc, Cor 1.6], it was proved that for a log del Pezzo
surface $Y$, the $Y^0$ is rationally connected and hence has
finite $\pi_1(Y^0)$ (see [Ca] and [KoMiMo]).

\par \vskip 1pc \noindent
{\bf Theorem 2.8} [Z7, Theorem 2].

\par \noindent
{\it Conjecture $2.4$ is true if one of the following occurs.}

\par \noindent
(1) $\dim V \le 2$.

\par \noindent
(2) {\it The Fano index $r(V) > \dim V -2$.}

\par \noindent
(3) $V$ has only isolated singularities and $r(V) = \dim V -2 = 1$.

\par \vskip 1pc \noindent
{\bf Theorem 2.9} [Z7, Theorem 2].

\par \noindent
{\it Let $V$ be a Fano variety of Fano index $r(V) > \dim V  - 2$ and
with at worst canonical singularities.
Then $\pi_1(V^0)$ is abelian of order $\le 9$.}

\par \vskip 1pc \noindent
{\bf Theorem 2.10} [Z7, Theorems 1 and 2].

\par \noindent
{\it Let $V$ be a Fano variety. Then $\pi_1(V^0) = (1)$ if one of the following
occurs:}

\par \noindent
(1) {\it The Fano index $r(V) > \dim V - 1$.}

\par \noindent
(2) {\it $\dim V = 3$ and $V$ has only Gorenstein isolated singularities.}

\par \vskip 1pc \noindent
{\bf (2.11).}
The following gives a concrete upper bound for $\pi_1(V^0)$ in certain case.
A relation $m(K_V+H) \sim 0$ in the theorem below occurs
when $V$ has Fano index 1 and Cartier
index $m$. It is conjectured that $m = 1, 2$.
To prove the theorem below, we show first
that there is a natural surjective map
$\pi_1(H^0) \rightarrow \pi_1(V^0)$ and also
use the fact that $H$ is Du Val K3 or Enriques.
Now the theorem follows from the results on $H$ in [KZ, Theorems 1 and 2].
For each of the three exceptional cases of $(p, c)$ below, we note
that there is a Du Val K3 or Enriques surface $Y$ with $\text{\rm Sing} \, Y = cA_{p-1}$
and $\pi_1(Y^0)$ infinity [ibid.].

\par \vskip 1pc \noindent
{\bf Theorem} [KZ, Theorem 3].

\par \noindent
{\it Let $p$ be a prime number. Let $V$ be a log terminal Fano $3$-fold with
a Cartier divisor $H$ such that $m(K_V+H) \sim 0$ (linear equivalence)
for $m = 1$ or $2$. Suppose that a member $H$ of $|H|$ is irreducible normal
and has $c$ singularities of type $A_{p-1}$ and no other singularities.}

\par \noindent
{\it Then $\pi_1(V^0)$ is soluble; and if $(p, c) \ne (2, 8), (2, 16), (3, 9)$,
then $|\pi_1(V^0)| \le 2p^k$ for some $k \le 4$.}

\par \vskip 1pc \noindent
{\bf Remark 2.12.}

\par \noindent
(1) In {\bf (2.4)} if we replace ''log terminal''
by ''log canonical'', then {\bf (2.4)} has counter-examples;
more precisely, if $V$ is a normal Fano surface with
at worst rational log canonical singularities, then
$\pi_1(V^0)$ contains a finite-index abelian subgroup of rank $k$
($k = 0, 2$) [Z8, Theorem 2.3].

\par \noindent
(2) In {\bf (2.9)}, the upper bound is optimum [MZ1, Lemma 6];
also ''canonical'' can not be replaced by ''log terminal''
[Z3, Appendix].

\par \vskip 1pc \noindent
{\bf (2.13).} In [Kj2], log del Pezzo surfaces with a unique
singularity are classified (including the existence part).
The classification of log del Pezzo surface of Cartier index $\le 2$
were announced in [AN]. In [Kj4], Kojima classified Picard number
1 log del Pezzo surfaces $Y$ of index 2, in a way different from [AN]:
there are exactly 18 types of $\text{\rm Sing} \, Y$ and the $\pi_1(Y^0) \le 8$;
the $\pi_1(Y^0) = (1)$ holds if and only if $Y$ contains the affine plane
as a Zariski open set.

\par \vskip 1pc \noindent
{\bf (2.14).} In [Ni3], the Picard number $\rho({\widetilde Y})$
of the minimal resolution ${\widetilde Y}$ of a log del Pezzo
surface $Y$ is bounded from above in terms of the maximum of
multiplicities of $Y$.

\head
Section 3. Normal algebraic surfaces $Y$ with Kodaira dimension $\kappa(Y^0) = 0$
\endhead

In this section we consider projective surfaces $Y$ with at worst
quotient singularities and Kodaira dimension $\kappa(Y^0) = 0$,
where $Y^0 = Y -$ \text{\rm Sing} \, $Y$.

\par \vskip 1pc \noindent
{\bf Theorem 3.1} ([Ka1, Theorem 2.2], [Mi, Ch II, 6.1.3]).

\par \noindent
{\it Let $Y$ be a projective surface with at worst quotient singularities.
Then the following are equivalent:}

\par \noindent
(1) {\it $Y$ is relatively minimal with Kodaira dimension $\kappa(Y^0) = 0$.}

\par \noindent
(2) {\it There is a positive integer $m$ such that
$mK_Y \sim 0$ (linear equivalence).}

\par \vskip 1pc \noindent
{\bf (3.2).} The smallest positive integer $m$ with $mK_Y \sim 0$
is called the {\it index} of $Y$ and denoted by $I = I(Y)$.

\par \vskip 1pc \noindent
{\bf Proposotion 3.3.} {\it Let $Y$ be a projective surface with at worst quotient
singularities and $I K_Y \sim 0$, where $I > 0$ is the index of $Y$.
Suppose that $Y$ is irrational. Then one of the following occurs.}

\par \noindent
(1) {\it $Y$ is a (smooth) abelian surface ($I = 1$) or a
hyperelliptic surface ($I = 2, 3, 4, 6$).}

\par \noindent
(2) {\it $Y$ has at worst Du Val singularities. The minimal resolution of $Y$
is either a $K3$ surface ($I = 1$) or an Enriques surface ($I = 2$).}

\par \vskip 1pc \noindent
{\it Proof.} In notation of {\bf (1.3)}, we have $I(K_{\widetilde Y} + D^*) \sim 0$.
So $\kappa({\widetilde Y}) \le 0$. If $\kappa({\widetilde Y}) = 0$,
Then {\bf (3.3)} follows from the classification of smooth surfaces.
If $\kappa({\widetilde Y}) = -\infty$, then ${\widetilde Y}$ is an irrational
ruled surface over a base curve of genus $\ge 1$.
However, $K_{\widetilde Y} + D^* = 0$ (numerically)
implies that $D^*$ contains some horizontal components
(this can be seen by going to a relative minimal model of ${\widetilde Y}$) which dominates
the base curve and hence is irrational, contradicting the fact that $D$
consists of rational curves only {\bf (1.3)}. So $\kappa({\widetilde Y}) = -\infty$
is impossible. This proves the proposition.

\par \vskip 1pc
In view of {\bf (3.3)}, to classify those $Y$
with at worst quotient singularities and $\kappa(Y^0) = 0$,
we need only to consider rational surfaces $Y$ with $mK_Y \sim 0$ for
some integer $m \ge 2$ (see {\bf (1.5)}). These are precisely
rational log Enriques surfaces {\bf (1.5)}.

\par \vskip 1pc \noindent
{\bf (3.4)}. Let $Y$ be a rational log Enriques surface.
Then the index $I = I(Y) \ge 2$.
Since $IK_Y \sim 0$,
there is a {\it canonical} ${\bold Z}/(I)$-Galois cover
$\pi: X = Spec \oplus_{j=0}^{I-1} {\Cal O}(-jK_Y) \rightarrow Y$
which is unramified over $(Y \setminus \{$non-Du Val singularities$\}) \supseteq Y^0$
and satisfies $K_X \sim 0$. Therefore, either $X$ is a (smooth) abelian surface
or a K3 surface possibly with some Du Val singularities.

\par
Recently, Suzuki [Su] has proved Morrison's cone conjecture for
rational log Enriques surfaces $Y$ : there is a finite rational polyhedral
cone which is a fundamental domain for the action of \text{\rm Aut}$(Y)$ on the rational
convex hull of its ample cone.

\par \vskip 1pc \noindent
{\bf Theorem 3.5.} {\it Suppose that $Y$ is a rational log Enriques surface
of index $I$ and that the canonical ${\bold Z}/(I)$-cover $X$ of $Y$
is an abelian surface. Then we have:}

\par \noindent
(1) {\it $I = 3$ or $5$. If $I = 3$, then $\text{\rm \text{\rm Sing} \,} Y$ consists of $9$ singularities of
type $(1/3)(1,1)$; if $I = 5$, then $\text{\rm \text{\rm Sing} \,} Y$ consists of $5$ singularities of
type $(1/5)(1, 2)$;} see [Re1] for notation.

\par \noindent
(2) {\it For each $I = 3, 5$, there is a unique
log Enriques surface $Y_I$ with $X_I$ abelain and $I(Y_I) = I$}.
To be precise, $Y_I = X_I/\langle g_I \rangle$, where
$X_3 = E_{\zeta_3} \times E_{\zeta_3}$ with
$E_{\zeta_3} = {\bold C}/({\bold Z} + {\bold Z} \zeta_3)$
an elliptic curve of period $\zeta_3 = \text{\rm exp}(2 \pi \sqrt{-1}/3)$,
$g_3 = \text{\rm diag}(\zeta_3, \zeta_3)$,
$X_5$ is the Jacobian surface of the genus-$2$ curve : $y^2 = x^5-1$,
and $g_5$ in $\text{\rm \text{\rm Aut}}(X_5)$ is induced by the curve automorphism :
$(x, y) \mapsto (\zeta_5x, y)$
(see [Bl, Example $1.2$], [Su, Proposition 1.2] and [Z4, Example $4.2$]).

\par \vskip 1pc \noindent
{\it Proof.} (1) is proved in [Z4, Theorem 4.1]. (2) is proved in [Bl, Su].

\par \vskip 1pc \noindent
{\bf Theorem 3.6.} {\it Let $Y$ be a log Enriques surface of index $I$.
Then $I \le 21$.}

\par \vskip 1pc \noindent
{\bf Remark 3.7.}

\par \noindent
(1) It is easy to see that the Euler function
$\varphi(I) \le 21$ and hence $I \le 66$ [Z4, Part I, Lemma 2.3].
In [Bl, Theorem C], it is proved that $I \le 21$.

\par \noindent
(2) Examples of $Y$ with prime $I(Y)$ are constructed in [Z4, Example 5.3-5.8;
Part II, Example 7.3] and [Bl, Example 4.1].

\par \vskip 1pc \noindent
{\bf (3.8).} A log Enriques surface $Y$ is {\it maximum} if any birational morphism
$Y \rightarrow Z$ to another log Enriques surface is an isomorphism.
By [Z4, Part II, Theorem 2.11'], for every log Enriques surface $Y$ of prime index,
there is a unique maximum log Enriques
surface $Y_{max}$ with $I(Y_{max}) = I(Y)$
and a birational morphism $Y_{max} \rightarrow Y$;
each singularity (if exists) of the canonical cover of $Y_{max}$ is of type $A_1$.

\par
The surface $Y(A_{19})$ in the assertion(3) below is not isomorphic to
the unique (modulo projective transformation) quartic $K3$ surface
in ${\bold P}^3$ with a Dynkin type $A_{19}$ singularity;
neither can $Y(D_{19})$ be embedded in ${\bold P}^3$ [KN].
The $Y(D_{19})$ is constructed in two different ways in
[Z4, Example 6.11: the $\overline{V}'$]
and [OZ1, Example 1], and $Y(A_{19})$ in [Z4, Example 3.2 : the $\overline{V}$]
and [OZ1, Example 2].
The uniqueness problem of $Y(D_{19})$ was initiated by
[Re2, Round 3, Example 6].

\par \vskip 1pc \noindent
{\bf Theorem} [OZ4, Corollary 4; OZ3, Corollary in \S 1; OZ1, Theorems 1 and 2].

\par \noindent
(1) {\it For each $I = 13, 17, 19$, there is a unique
maximum log Enriques surface $Y$ with $I(Y) = I$.}

\par \noindent
(2) {\it All maximum log Enriques surfaces of index $11$ form
a family of dimension $1$ and are all given in} [OZ3, \S 1, Corollary].

\par \noindent
(3) {\it For each $D$ in $\{D_{19}, A_{19}\}$, there is a unique
rational log Enriques surface $Y(D)$ whose canonical cover has a singularity
of Dynkin type $D$. The index $I(Y)$ equals $3$ (resp. $2$) when $Y$ equals
$Y(D_{19})$ (resp. $Y(A_{19})$).}

\par \vskip 1pc
Next we will investigate the behaviour of $\pi_1(Y^0)$ and propose
a conjecture {\bf (3.12)} generalizing the one in [CKO]. Note that
${\bold Z}/(I)$ is the image of $\pi_1(Y^0)$ by a homomorphism.

\par \vskip 1pc \noindent
{\bf Theorem 3.9} [Z4, Part II, Theorem 2.11', Cor 1]. {\it Let $Y$ be
a maximum log Enriques surface of odd prime index $I$.
Let $X \rightarrow Y$ be the
canonical ${\bold Z}/(I)$-cover. Then we have:}

\par \noindent
(1) {\it $X$ has at worst type $A_1$ singularities and
$\#(\text{\rm Sing} \, X) \le 6$.}

\par \noindent
(2) $\pi_1(Y^0) = {\bold Z}/(I)$.

\par \vskip 1pc \noindent
{\it Proof.} (1) is proved in [Z4, Part II, Cor. 1].
For (2), we have only to show that $\pi_1(X^0) = (1)$
because the inverse of $Y^0$ via the canonical map $X \rightarrow Y$
(unramified over $Y^0$)
is $X^0$ with a few smooth points removed (the removal of smooth
points in a complex surface does not change $\pi_1$).
Since $\#(\text{\rm Sing} \, X) \le 6$ by (1), we have $\pi_1(X^0) = (1)$ [KZ, Theorem 1].
This proves the theorem.

\par \vskip 1pc \noindent
{\bf Theorem 3.10} [SZ, Proposition 4.1, Theorem 4.3, Cor 4.4]

\par \noindent
(1) {\it Suppose that "$\pi_1^{alg}(X^0) = (1) \Rightarrow \pi_1(X^0) = (1)$
for all Du Val $K3$ surfaces $X$" holds ($X$ satisfies this condition
if $(*)$ in $(2)$ holds for $X$).
Let $Y$ be a log Enriques surface of index $I$.
Then either $\pi_1(Y^0)$ is finite or there is a finite morphism
$Z \rightarrow Y$ from an abelian surface which is unramified over $Y^0$.
In particular, $\pi_1(Y^0)$ contains a finite-index abelian subgroup
of rank $k$ ($k = 0, 4$).}

\par \noindent
(2) {\it Let $X \rightarrow Y$ be the canonical cover of a log Enriques surface
of index $I$,
let ${\widetilde X} \rightarrow X$ be the minimal resolution
and $D = \sum D_i$ the exceptional divisor.
Then $\pi_1(Y^0) = {\bold Z}/(I)$ if
the lattice $\Gamma = {\bold Z}[\cup D_i]$ is primitive in
$H^2({\widetilde X}, {\bold Z})$ and satisfies:}

\par \noindent
(*) {\it $r = rank(\Gamma) \le 18$ and the discriminant group
$\Gamma^{\vee}/\Gamma$
is generated by $k$ elements with $k \le min\{r, 20-r\}$.}

\par \vskip 1pc \noindent
{\it Proof.} Let $X \rightarrow Y$ be the canonical ${\bold Z}/(I)$-cover.
Now the conclusion in {\bf (3.10)} (1) with $Y$ replaced by $X$ holds
by [SZ, Proposition 4.1]. So (1) is true for $X \rightarrow Y$
is unramified over $Y^0$. For (2), [SZ, Theorem 4.3] implies
$\pi_1(X^0) = (1)$. So (2) is true. This proves the theorem.

\par \vskip 1pc \noindent
There is a concrete upper bound of $\pi_1(Y^0)$ for certain $Y$.

\par \vskip 1pc \noindent
{\bf Theorem 3.11} [GZ3, Theorem 1].
{\it Let $Y$ be a rational log Enriques surface of index $2$.
Assume that $Y$ has no Du Val singularities. Then $\pi_1(Y^0)$
is a soluble group of order $n_1 n_2$ with $n_i \le 16$.}

\par \vskip 1pc
The results {\bf (3.9)} $\sim$ {\bf (3.11)} support the following
which is just the conjecture in [CKO]
when $Y$ is a Du Val $K3$ surface, i.e., when $I(Y) = 1$ and $q(Y) = 0$.

\par \vskip 1pc \noindent
{\bf Conjecture 3.12.} {\it Let $Y$ be a log Enriques surface.
Then the universal cover $U$ of $Y^0$ is a big open set (= the
complement of a discrete subset) of either a Du Val $K3$ surface
or of ${\bold C}^2$; in the latter case, $U \rightarrow Y^0$
factors through a finite etale cover $Z^0 \rightarrow Y^0$, where
$Z^0$ is a big open set of an abelian surface $Z$.}

\par \vskip 1pc \noindent
{\bf (3.13).} Since the canonical cover $X \rightarrow Y$ of
a log Enriques surface is unramified over $Y^0$, we
have $\pi_1(Y^0)/\pi_1(X^0) = {\bold Z}/(I)$. So
{\bf (3.12)} is, in most cases, reduced to the problem of $\pi_1(X^0)$
for a Du Val K3 surface $X$ (see {\bf (3.5)}). We have the following:

\par \vskip 1pc \noindent
{\bf Theorem} [KZ, Theorems 1 and 2]. {\it Conjecture $3.12$
is true if $Y$ is a Du Val $K3$ or Enriques surface and has
several singularities of type $A_{p-1}$ and no other singularities;
here $p$ is a prime number.}

\par \vskip 1pc \noindent
{\bf (3.14).} The following results contributes towards an
answer in affirmative to {\bf (3.12)}.
These are just applications of [CKO, Theorems A and B] to
the canonical cover $X$ of $Y$. Also the upper bound
$\#D \le 15$ is an optimum condition for $\pi_1(X^0)$ to be finite
by considering Kummer surfaces.

\par \vskip 1pc \noindent
{\bf Theorem.} (1) {\it Let $Y$ be a log Enriques surface with an elliptic
fibration. Then either $\pi_1(Y^0)$ is finite or there is a finite
cover of $Z \rightarrow Y$ from an abelian surface which is unramified over $Y^0$.}

\par \noindent
(2) {\it Let $Y$ be a log Enriques surface, $X \rightarrow Y$ the canonical
cover and ${\widetilde X} \rightarrow X$ the minimal resolution
with $D$ the exceptional divisor. Suppose that $\#D \le 15$.
Then $\pi_1(Y^0)$ is finite.}

\par \vskip 1pc \noindent
{\bf (3.15).} In [Oh], pairs $(S, \Delta)$ of normal surface $S$
and a ${\bold Q}$-divisor $\Delta$ satisfying $K_S + \Delta \equiv 0$
(numerically) are considered. These pairs appear naturally as degenerate
fibres in log degeneration; for many interesting cases, he completed the
classification of these pairs.

\par \vskip 1pc \noindent
{\bf (3.16).} In [Kj1, Theorem 0.1], strongly minimal smooth affine
surface $S$ with $\kappa(S) = 0$ is classified and its invariants
are classified ({\it strongly minimal} means almost minimal and having no
exceptional curve of the second kind [Mi, Ch II, {\bf (4.9)}]).
In particular, the minimal $m > 0$ with log pluri-genus $P_m(S) > 0$,
the log irregularity $q(S)$ and the Euler number $e(S)$ satisfy
the following ($\pi_1(S)$ is also calculated there, which is generated
by at most two elements):
$$m | 6, \,\,\, q(S) \in \{0, 1, 2\}, \,\,\, e(S) \in \{0, 1, 2, 3, 4\}.$$

\par \vskip 1pc \noindent
{\bf (3.17).} In [Fu, \S 8], all {\bf Q}-homology planes of Kodaira dimension 0
are classified. It was also proved there that there is no
{\bf Z}-homology plane $S$ of Kodaira dimension $\kappa(S) = 0$. The
paper [Fu] is very important and also essentially used in [Kj1].

\par \vskip 1pc \noindent
{\bf (3.18).} Iitaka [I2] conjectured that an affine normal variety
$S$ is isomorphic to $({\bold C}^*)^n$ if and only if $\kappa(S) = 0$
and $q(S) = \dim S$. In the same paper, he himself proved it when $\dim S = 2$.

\par
According to [I1], a (possibly open) surface $S$ is {\it logarithmic $K3$} if
the logarithmic invariants satisfy : $q(S) = 0$, $p_g(S) = 1$, $\kappa(S) = 0$.
In [I1] log K3 surfaces were classified. In [Z1], one defines
the {\it Iitaka surface} as
a pair $(V, A+N)$ of smooth projective rational surface $V$
and reduced divisor $A+N$ with $A+K_V \sim 0$ and $N$ contractible
to Du Val singularities, and the classification of such pairs were done there.

\par \vskip 1pc \noindent
{\bf (3.19).} In [Kj3], Kojima studies complements $S$ of reduced plane curves
with $\kappa(S) = 0$; in particular
he proves that the logarithmic geometric genus $p_g(S) = 1$.

\head
Section 4. Normal algebraic surfaces $Y$ with Kodaira dimension
$\kappa(Y^0) = 1, 2$
\endhead

In this section we consider projective surface $Y$
with at worst quotient singularities and $\kappa(Y^0) = 1, 2$.

\par \vskip 1pc \noindent
{\bf (4.1).} We first consider the case $\kappa(Y^0) = 1$. The
following is a consequence of [Ka1, Theorem 2.3] or [Mi, Ch II,
Theorem 6.1.4]. Indeed, in our case, the boundary divisor $D^*$
is fractional and contains no effective integral divisor {\bf
(1.3)}.

\par \vskip 1pc \noindent
{\bf Theorem.} {\it Let $Y$ be a projective surface with at worst
quotient singularities.
Suppose that $Y$ is relatively minimal and $\kappa(Y^0) = 1$.}

\par \noindent
{\it Then there is a positive integer $m$ such that
$mK_Y$ is Cartier and the linear system $|mK_Y|$ is composed with an irreducible
pencil  $\Lambda$ without base points. Each general member of $\Lambda$
is a smooth elliptic curve. So there is an elliptic fibration $Y \rightarrow B$.}

\par \vskip 1pc \noindent
{\bf (4.2).} Let $Y$ be a projective surface with at worst quotient singularities.
Suppose that
$Y$ is relatively minimal and $\kappa(Y^0) = 2$. Then $K_Y$ is nef and big.
By [KMM, Theorem 3-1-1], $|mK_Y|$ is base point free for $m$ sufficiently
divisible,
and hence defines a birational morphism $\varphi : Y \rightarrow Z$.
This $\varphi$ is nothing but the contraction of all curves
on $Y$ having zero intersection with $K_Y$.
Then $Z$ has at worst log canonical singularities
([Ka1, Theorem 2.9], [Mi, Ch II, Theorem 4.12]).
Denote by $LC$ the set of points on $Z$ which is log canonical but not log terminal
(i.e., not of quotient
singularity). Then we have the following Miyaoka-Yau type inequality proved by
[Kb1, 2] and [KNS].

\par \vskip 1pc \noindent
{\bf  Theorem 4.3.} {\it Let $Y$ be a projective surface
with at worst quotient singularities. Suppose that
$Y$ is relatively minimal and $\kappa(Y^0) = 2$.
Then we have the following, where $P$ runs over all quotient singularities
of $Z$ and $G_P$ is the local fundamental group at $P$}
$$K_Y^2 \le 3\{e(Z - LC) - \sum_P (1 - \frac{1}{|G_P|})\}.$$

\par \vskip 1pc \noindent
{\bf (4.4).} For smooth and minimal surfaces $X$ of general type, we
have Noether inequalities
$$p_g(X) \le (1/2)c_1(X)^2 + 2, \,\,\,\, c_1(X)^2 \ge \frac{1}{5}c_2(X) - \frac{36}{5}.$$

For singular surfaces, we have:

\par \vskip 1pc \noindent
{\bf Theorem} [TZ, Theorems 1.3 and 2.10].
\par \noindent
{\it Let $Y$ be a projective surface with at worst quotient
singularities and $\kappa(Y^0) = 2$. Then the logarithmic
geometric genus $p_g(Y^0)$ satisfies the optimum upper bound:}
$$p_g(Y^0) < K_Y^2 + 3.$$

\par \vskip 1pc \noindent
{\bf (4.5).}  For smooth projective surface $X$ of general type,
the famous Miyaoka-Yau inequality asserts that
$c_1(X)^2 \le 3c_2(X)$.

\par
Consider log surface $(V, D)$ with
$V$ a smooth projective surface and $D$ a reduced divisor with
simply normal crossings.
Set $\overline{c}_1^2 = (K_V+D)^2$ and $\overline{c}_2 = c_2(V) - e(D)$.
Sakai [Sa] proved that
$$\overline{c}_1^2 \le 3\overline{c}_2 $$
provided that
$D$ is semi-stable and $\kappa(V \setminus D) = 2$.
The following is a lower bound of $\overline{c}_1^2 $ in terms
of $\overline{c}_2 $. These two inequalities
together give effective restrictions on the region for non-complete
algebraic surfaces $V \setminus D$ of general type to exist.
In the following, $(V, D)$ is {\it minimal} if
$K_V+D$ is nef and there is no $(-1)$-curve $E$ with $E . (K_V+D) = 0$.

\par \vskip 1pc \noindent
{\bf Theorem} [Z5, Cor. to Theorem C, Theorem D]. {\it Let $(V, D)$ be a log surface
with $D \ne 0$.
Assume that $(V, D)$ is minimal and $\kappa(V \setminus D) = 2$.
Assume further that $\kappa(V) \ge 0$. Then we have}

\par \noindent
(1)
$$\overline{c}_1^2  \ge \frac{1}{15} \overline{c}_2  - \frac{8}{5}.$$

\par \noindent
(2) {\it Suppose that $p_g(V \setminus D) \ge 3$ and $|K_V+D|$ is
not composed with a pencil. Then}
$$\overline{c}_1^2  \ge \frac{1}{9} \overline{c}_2 - 2.$$

\head
Section 5. Automorphisms of algebraic surfaces - smooth surface case
\endhead

\par \noindent
{\bf (5.1).} We mention some background of $\text{\rm Aut}(X)$ where $X$
is a smooth projective rational surface.
$\text{\rm Aut}(X)$ had been studied by S. Kantor more than
one hundred years ago [Kt].
It was continued by Segre, Manin, Iskovskih, Gizatullin
and many others [Se], [Ma1, 2], [Is], [Gi].
See also [Ho1], [Ho2].
In [DO], the group of automorphisms of any general del Pezzo surface is described
and it turns out that its discrete part is equal to the
kernel of the Cremona representation on the moduli space of $n$ points in ${\bold P}^2$.
Very recently, de Fernex [dF] constructed all the Cremona transformations of ${\bold P}^2$
of prime order, where he employed the methods different from those used
by Dolgachev and Zhang in [ZD].

\par \noindent
In [ZD], minimal pairs $(X, G)$ with prime order $p = |G|$ was considered.
In particular, using the recent Mori theory,
it was shown there that if the $G$-invariant sublattice
of $\text{\rm \text{\rm Pic}  \, } X$ has rank $1$ then $p \le 5$
unless $X = {\bold P}^2$;
the short and precise classification of these pairs, modulo equivariant
isomorphism, was also given there.

\par \vskip1pc \noindent
Generic Enriques surfaces have infinitely many automorphisms.
Those Enriques surfaces with finite automorphisms have been
classified by S. Kond$\bar{\roman o}$ ([Kon4], see also [Ni2]);
there are seven families of such Enriques surfaces. It follows
that K3 covers of Enriques surfaces all have infinite
automorphism groups.

\par \vskip1pc \noindent
For K3 surfaces, much progress on their automorphism groups has
been done by Kond$\bar{\roman o}$, Keum, Dolgachev, Oguiso, and
Zhang ([Kon1, 2, 3], [Ke1, 3], [KK], [DK], [OZ1, 3, 4, 5]).

\par \vskip 1pc \noindent
It is known that minimal surfaces of general type has only finite
automorphism groups.  It was Xiao [Xi1] who gave a proof of the
existence of a bound for the order of the automophism group,
which is linear in the Euler number of the surface. For curves
$C$ of genus $\ge 2$, a classical theorem of Hurwitz gave a sharp
bound $|\text{\rm Aut}(C)|\le 84(g(C)-1)=-42e(C)$.

\par \vskip 1pc \noindent
{\bf Theorem} [Xi2, Theorem 2]. {\it Let $X$ be a minimal surface
of general type. Then $|\text{\rm Aut}(X)|\le (42K_X)^2$, with equality if
and only if $X\cong (C\times C)/N$, where $C$ is a curve with
$|\text{\rm Aut}(C)|=84(g(C)-1)$, $N$ a normal subgroup of $\text{\rm Aut}(C\times C)$
acting freely on $C\times C$ and preserving the two projections of
$C\times C$.}

\par \vskip 1pc \noindent
{\bf (5.2).} In [ZD], pairs $(X, G)$ of a smooth {\it rational} projective
surface $X$ and a finite group of automorphisms are considered.
A pair is {\it minimal} if every $G$-equivariant birational morphism to another pair
$(Z, G)$ is an isomorphism.

\par \vskip 1pc \noindent
{\bf Theorem} [ZD, Therorem 1]. {\it Let $(X, \mu_p)$
be a minimal pair with $p$ a prime number.}

\par \noindent
(1) {\it If the invariant sub-lattice $(\text{\rm Pic}  \,  X)^{\mu_p}$ has rank at least $2$,
then $X$ is a Hirzebruch surface and the pair $(X, \mu_p)$ is birationally
equivalent to a pair $({\bold P}^2, \mu_p)$.}

\par \noindent
(2) {\it If $(\text{\rm Pic}  \,  X)^{\mu_p}$ has rank $1$, then the pairs are classified
in $[ZD, Theorem 1]$; in particular, we have $p \le 5$ unless $X = {\bold P}^2$.}

\par \vskip 1pc \noindent
{\bf (5.3).} Let $X$ be a K3 surface. The following are well
known (cf. [BPV]).

\par \noindent
(1) $H^{2,0}(X)={\bold C}\omega_X$,

\noindent
where $\omega_X$ is a nowhere vanishing global
holomorphic 2-form on $X$.

\par \noindent
(2) $H^2(X,{\bold Z})$ is an even unimodular lattice of signature
(3,19) with the cup product, so we have an isomorphism
$$H^2(X,{\bold Z})\cong U\oplus U\oplus U\oplus E_8(-1) \oplus
E_8(-1),$$ where $U$ (resp. $E_8$) is the even unimodular lattice
of signature (1,1) (resp. (8,0)).

\par \noindent
(3) $\text{\rm Pic}  \, (X)$ is isomorphic to the N\'eron-Severi group $NS(X)$, and
hence can be viewed as a sublattice of $H^2(X,{\bold Z})$. The
rank of $\text{\rm Pic}  \, (X)$, called the Picard number of $X$, is denoted by
$\rho(X)$. This number can take the value 0, 1, ..., 20. The
lattice $\text{\rm Pic}  \, (X)$ is hyperbolic(=Lorentzian) if $X$ is projective,
and is semi-negative definite or negative definite if $X$ is not
projective.

\par \noindent
(4) Recall that all K3 surfaces are K\"{a}hler [Siu], so Hodge
decomposition holds for them.

\par \noindent
(5) Let $$C(X)\subset H^{1,1}(X, {\bold R}):= H^{1,1}(X)\cap
H^2(X, {\bold R})$$ denote the K\"{a}hler cone of $X$, the set of
all classes of symplectic forms of K\"{a}hler-Einstein metrics on
$X$. In K3 surface case, $C(X)$ can be numerically characterized
as follows:
\medskip
$C(X)=\{\omega\in H^{1,1}(X, {\bold R}) :
\langle\omega,\omega\rangle>0, \langle\omega,R\rangle>0$ for all
smooth rational curves $R\}$

\par \vskip 1pc \noindent
For a compact K\"{a}hler manifold, Nakai-Moishezon type
criterion, i.e. the characterization of the K\"{a}hler cone, is
highly non-trivial (see [DP]).

\par \vskip 1pc \noindent
(6) Let $r$ be an element of $\text{\rm Pic}  \, (X)$ with $\langle r,
r\rangle=-2d$ $(d>0)$, and $ \langle r,\text{\rm Pic}  \, (X)\rangle\subset
d{\bold Z}.$ Then $$x\to x+\langle x,r\rangle r/d$$ defines an
isometry of $\text{\rm Pic}  \, (X)$, called a $(-2d)$-reflection.

\par
Let $W(\text{\rm Pic}(X))$ (resp.$W(\text{\rm Pic}(X))^{(2)}$)
be the subgroup of the
orthogonal group $O(\text{\rm Pic}(X))$ generated by all reflections (resp.
all $(-2)$-reflections). These are normal subgroups of $O(\text{\rm Pic}  \, (X))$
and, by linearity, acts naturally on $H^{1,1}(X, {\bold R})$. The
set $$\{\omega\in H^{1,1}(X, {\bold R}) :
\langle\omega,\omega\rangle>0\}$$ has two components, each a cone
over a 19-dimensional hyperbolic manifold with constant
curvature. $C(X)$ is contained in one of the two components, and
the action of $W(\text{\rm Pic}  \, (X))^{(2)}$ on this component has $C(X)$ as
its fundamental domain.

\par \vskip 1pc \noindent
(7) If $X$ is projective, the ample cone $$D(X):= C(X)\cap
\text{\rm Pic}  \, (X)\otimes{\bold R}$$ is non-empty and can be numerically
characterized as
\medskip
$D(X)=\{\omega\in \text{\rm Pic}  \, (X)\otimes{\bold R} :
\langle\omega,\omega\rangle>0, \langle\omega,R\rangle>0$ for all
smooth rational curves $R\}$.
\medskip \noindent
The group $W(\text{\rm Pic}  \, (X))^{(2)}$ acts on the component $P_+(X)$ of
$$\{\omega\in \text{\rm Pic}  \, (X)\otimes{\bold R} :
\langle\omega,\omega\rangle>0\}$$ containing $D(X)$, and has
$D(X)$ as its fundamental domain.

Note that $O(\text{\rm Pic}  \, (X))$ acts on $P_+(X)$ and is a semi-direct
product of the normal subgroup $W(\text{\rm Pic}  \, (X))^{(2)}$ and the symmetry
group $SymD(X)$ of the cone $D(X)$, i.e.
$$O(\text{\rm Pic}  \, (X))/W(\text{\rm Pic}  \, (X))^{(2)}\cong SymD(X).$$

\par \vskip 1pc \noindent
{\bf (5.4).} The Torelli theorem asserts that a K3 surface is
determined up to isomorphism by its Hodge structure. More
precisely we have:

\par \vskip 1pc \noindent
{\bf Theorem} ([PSS],[BR]). {\it Let $X$ and $Y$ be K3 surfaces,
and let $$\phi : H^2(X,{\bold Z})\to H^2(Y,{\bold Z})$$ be an
isometry. Extend $\phi$ to $H^2(X,{\bold C})$ or to $H^2(X,{\bold
R})$ by tensoring with ${\bold C}$ or ${\bold R}$. Then :}

\par \noindent
(1) {\it If $\phi$ sends $H^{2,0}(X)$ to $H^{2,0}(Y)$, then $X$
and $Y$ are isomorphic.}

\par \noindent
(2) {\it If $\phi$ also sends $C(X)$ to $C(Y)$, then $\phi=f^*$
for a unique isomorphism $f:Y\to X$.}

\par \vskip 1pc \noindent
{\bf (5.5).} Let $X$ be a projective K3 surface. Torelli theorem
shows that there is a map
$$\text{\rm Aut}(X) \to O(\text{\rm Pic}  \, (X))/W(\text{\rm Pic}  \, (X))^{(2)}\cong SymD(X)$$
which has finite kernel and cofinite image. So in practice, if we
want to describe $\text{\rm Aut}(X)$, the main step is to calculate
$O(\text{\rm Pic}  \, (X))/W(\text{\rm Pic}  \, (X))^{(2)}$. This is in general a highly
nontrivial arithmetic problem, if the group is infinite. There
are 3 cases:

\par \vskip 1pc \noindent
(1) $W(\text{\rm Pic}  \, (X))^{(2)}$ is of finite index in $O(\text{\rm Pic}  \, (X))$.

\par \vskip 1pc \noindent
(2) $W(\text{\rm Pic}  \, (X))^{(2)}$ is of infinite index in $O(\text{\rm Pic}  \, (X))$, but
$W(\text{\rm Pic}  \, (X))$ is of finite index in $O(\text{\rm Pic}  \, (X))$. (In this case, we
call $\text{\rm Pic}  \, (X)$ reflective.)

\par \vskip 1pc \noindent
(3) $W(\text{\rm Pic}  \, (X))$ is of infinite index in $O(\text{\rm Pic}  \, (X))$, i.e. $\text{\rm Pic}  \, (X)$ is not reflective.

\par \vskip 1pc \noindent
{\bf Remark.} The case (1) occurs if and only if $\text{\rm
Aut}(X)$ is finite. If $\rho(X) \ge 3$, this occurs if and only
if $X$ contains at least one but finitely many smooth rational
curves. If $\rho(X)=2$, this occurs if and only if $X$ contains a
smooth rational curve or an irreducible curve of arithmetic genus
1, if and only if $\text{\rm Pic}  \, (X)$ represents $-2$ or 0
[PSS]. Nikulin [Ni1, Ni4] and Vinberg classified all such
lattices of rank $\ge 3$ belonging to the case (1). It follows
from the classification that every algebraic Kummer surface has an
infinite automorphism group (cf. [Ke2]).

\par
The classification of the N\'eron-Severi lattice is also utilized
in [Og1 - Og3], where he has proved the density of the jumping
loci of the Picard number of a hyperk\"ahler manifold under small
1-dimensional deformation, where he reveals the structure of
hierarchy among all the narrow Mordell-Weil lattices of Jacobian
K3 surfaces.

\par \vskip 1pc \noindent
{\bf (5.6).} For finite groups which can act on a K3 surface, the
following results are given by S. Mukai and S. Kondo.

\par \vskip 1pc \noindent
{\bf Theorem} [Mu],[Kon2]. {\it Let $X$ be a $K3$ surface and let
$G$ be a finite symplectic subgroup of $\text{\rm Aut}(X)$, i.e. $G$ acts
trivially on $H^{2,0}(X)$. Then $G$ is isomorphic to a subgroup
of the Mathieu group $M_{23}$, which has at least five orbits on
a set $\Omega$ of $24$ elements. In particular, $|G|\le 960$.}

\par \vskip 1pc \noindent
{\bf Theorem} [Kon3]. {\it The maximum order among all finite
groups which can act on a $K3$ surface is $3840$, and is uniquely
realized by the group $({\bold Z}/2{\bold Z})^4\cdot A_5\cdot
{\bold Z}/4{\bold Z}$ acting on the Kummer surface
$Km(E_{\sqrt{-1}} \times E_{\sqrt{-1}})$, where $E_{\sqrt{-1}}$
is the elliptic curve with $\sqrt{-1}$ as its fundamental period.}

\par \vskip 1pc \noindent
Some projective K3 surfaces, including all algebraic Kummer
surfaces and K3 covers of Enriques surfaces, have infinite
automorphism groups. Given a projective K3 surface $X$ with
$\text{\rm Aut}(X)$ infinite, it is an interesting problem to determine a
set of geometric generators of $\text{\rm Aut}(X)$. This problem has been
settled for certain classes of K3 surfaces. These results are
given in {\bf (5.7) - (5.10)} below.

\par \vskip 1pc \noindent
{\bf (5.7).} {\it Two most algebraic $K3$ surfaces}

\par
Vinberg[Vin] calculated $\text{\rm Aut}(X)$ for two K3 surface
with transcendental lattice
$$T(X) =
\pmatrix 2 & 1 \\ 1 & 2 \endpmatrix, \,\,\,\,
\pmatrix 2 & 0 \\ 0 & 2
\endpmatrix $$
respectively. In both cases, the full reflection group $W(\text{\rm Pic}  \, (X))$
is of finite index in $O(\text{\rm Pic}  \, (X))$.

\par \vskip 1pc \noindent
{\bf (5.8).}  {\it generic Jacobian Kummer surfaces}

\par
Let $C$ be a smooth curve of genus 2. The Jacobian variety $J(C)$
of $C$ is an abelian surface with a natural involution $\tau$ and
the quotient variety $J(C)/\tau$ has 16 singularities of type
$A_1$. This surface can be embedded as a quartic surface $F$ in
${\bold P}^3$ with 16 nodes. The minimal resolution $X$ of
$J(C)/\tau$ is called the {\it Jacobian Kummer surface}
associated with $C$. We call $X$ {\it generic} if the
N\'eron-Severi group of $J(C)$ is generated by the class of $C$.
For $X$ generic, the transcendental lattice $T(X)$ can be
computed as follows:
$$T(X) = U(2)\oplus U(2)\oplus <-4>.$$

\par \noindent
Note that $\text{\rm Aut}(X)$ is isomorphic to the birational automorphism
group $Bir(F)$ of the singular quartic surface $F$. At the last
century it was known that $X$ has many involutions, that is, {\it
sixteen translations} induced by those of $J(C)$ by a 2-torsion
point, {\it sixteen projections} of $F$ from a node, {\it sixteen
correlations} by means of the tangent plane collinear to a trope,
and a {\it switch} defined by the dual map of $F$. In 1900,
Hutchinson found another 60 involutions associated with G\"opel
tetrads. Since Hutchinson, for generic $X$ no other automorphism
had been provided until new 192 automorphisms were given in [Ke1].

\par \vskip 1pc \noindent
{\bf Theorem} [Ke1]. {\it For a generic Jacobian Kummer surface,
there are $192$ new automorphisms of infinite order which are not
generated by classical involutions.}

\par \vskip 1pc \noindent
{\bf Theorem} [Kon1]. {\it The automorphism group of a generic
Jacobian Kummer surface is generated by the classical involutions and
the $192$ new automorphisms.}

\par \vskip 1pc \noindent
{\bf Theorem} [Ke3]. {\it For $F$ generic, all birational
automorphisms of $F$ are induced by Cremona transformations of
${\bold P}^3$.}

\par \vskip 1pc \noindent
{\bf (5.9).} {\it Kummer surfaces associated with the product of
two elliptic curves}

\par \noindent
The following four cases were considered. In each case, a set of
generators of $\text{\rm Aut}(X)$ is given in [KK].

\medskip\noindent
{\it Case} I. $X = Km(E \times F)$ where $E$ and $F$ are
non-isogenus generic elliptic curves.

\smallskip\noindent
{\it Case} II. $X = Km(E \times E)$ where $E$ is an elliptic
curve without complex multiplications.

\smallskip\noindent
{\it Case} III. $X = Km(E_{\omega} \times E_{\omega})$ where
$\omega$ is a 3rd root of unity and $E_{\tau}$ is the elliptic
curve with $\tau$ as its fundamental period.

\smallskip\noindent
{\it Case} IV. $X = Km(E_{\sqrt{-1}} \times E_{\sqrt{-1}})$.

\par \vskip 1pc \noindent
The transcendental lattice $T(X)$ can be computed as follows:
$$T(X) = \cases U(2)\oplus U(2),&\text{(Case I)}; \\
U(2) \oplus <4>,&\text{(Case II)}; \\
\pmatrix 4 & 2 \\ 2 & 4 \endpmatrix,&\text{(Case III)}; \\
\pmatrix 4 & 0 \\ 0 & 4 \endpmatrix,&\text{(Case IV)}.\endcases $$

\par \vskip 1pc \noindent
{\bf Remark.} In Case I the group $W(\text{\rm Pic}  \, (X))$ is of finite index
in $O(\text{\rm Pic}  \, (X))$ and in other cases not.

\par \vskip 1pc \noindent
{\bf (5.10).} {\it Quartic Hessian surfaces}

Let $S:F(x_0,x_1,x_2,x_3) = 0$ be a nonsingular cubic surface in
${\bold P}^3$. Its Hessian surface is a quartic surface defined
by the determinant of the matrix of second order partial
derivatives of the polynomial $F$. When $F$ is general enough,
the quartic $H$ is irreducible and has 10 nodes. It contains also
10 lines which are the intersection lines of five planes in
general linear position. The union of these five planes is
classically known as the Sylvester pentahedron of $S$. The
equation of $S$ can be written as the sum of cubes of some linear
forms defining the five planes. A nonsingular model of $H$ is a
K3 surface $\tilde H$. Its Picard number $\rho$ satisfies the
inequality $\rho \ge 16$. Note that $\text{\rm Aut}(\tilde H)\cong Bir(H)$.
In [DK] an explicit description of the group $Bir(H)$ is given
when $S$ is general enough so that $\rho = 16$. In this general
case, the transcendental lattice $T(\tilde H)$ can be computed as
follows:
$$T(\tilde H) = U\oplus U(2)\oplus A_2(-2).$$
 Although $H$, in general, does not have any non-trivial
automorphisms (because $S$ does not), the group $Bir(H)\cong
\text{\rm Aut}(\tilde H)$ is infinite.  It is generated by the automorphisms
defined by projections from the nodes of $H$, a birational
involution which interchanges the nodes and the lines, and the
inversion automorphisms of some elliptic pencils on $\tilde H$.

This can be compared with the known structure of the group of
automorphisms of the Jacobian Kummer surface {\bf (5.8)}. Indeed,
the latter surface is birationally isomorphic to the Hessian $H$
of a cubic surface [Hu] but the Picard number of $\tilde H$ is
equal to 17 instead of 16.

\par \vskip 1pc \noindent
{\bf (5.11).} Let us explain the method for computing the
automorphism group of an algebraic K3 surface, which was first
employed by S. Kond$\bar{\roman o}$ for generic Jacobian Kummer
surface case {\bf (5.8)}[Kon1]. The two cases {\bf (5.9)} and {\bf
(5.10)} use this method, and even the  first case {\bf (5.7)} can
also be calculated by the same method (See [Bor2]).

Let $X$ be an algebraic K3 surface with large Picard number, say
$\rho(X) \ge 3$. Suppose that $\text{\rm Aut}(X)$ is infinite.
Then the ample cone $D(X)$ is not a (finite) polyhedral cone,
i.e. has infinitely many faces. Hence, it is difficult to
describe $D(X)$ explicitly. Assume that one can find

\par \vskip 1pc
$\bullet$ a polyhedral cone $D'$ in $D(X)$,

\par \vskip 1pc
$\bullet$ a set of automorphisms $\{g_\alpha\}$ of $X$ whose
action on $D(X)$ has $D'$ as a fundamental domain.

\par \vskip 1pc \noindent
Then, by {\bf (5.5)}, one can conclude that the automorphisms
$\{g_\alpha\}$ generate the whole group $\text{\rm Aut}(X)$, up to finite
groups. In addition to $\{g_\alpha\}$, some symmetries of $D'$
(not all elements of $SymD'$ in general) may realize as
automorphisms of $X$ and some projectively linear automorphisms,
if any, generate the kernel of the map in {\bf (5.5)}.

\par \vskip 1pc \noindent
{\bf Remark.} In the above, $g_\alpha$ corresponds to a face of
$D'$ orthogonal to a vector $\alpha$, and acts on $D(X)$ like a
reflection, i.e. sends one of the half-spaces defined by $\alpha$
to the other half-space defined by $\alpha$ or to one of the two
half-spaces corresponding to $g_\alpha^{-1}$. The second case
actually occurs in generic Jacobian Kummer surface case {\bf
(5.8)}.

\par \vskip 1pc \noindent
{\bf (5.12).} To find such a polyhedral cone $D'\subset D(X)$,
Kond$\bar{\roman o}$ used the known structure of the orthogonal
group of the even unimodular lattice $II_{1,25}$ of signature
(1,25). (Such a lattice is unique up to isomorphism and is
isomorphic to $\Lambda\oplus U$, where $\Lambda$ is the Leech
lattice, i.e. the even unimodular negative definite lattice of
rank 24 which contains no vectors of norm $-2$.) To be more
precise, the following steps lead to the calculation of $\text{\rm Aut}(X)$.

\par \vskip 1pc \noindent
{\bf Step 1.} Compute $\text{\rm Pic}  \, (X)$.

\par \vskip 1pc \noindent
{\bf Step 2.} Embed $\text{\rm Pic}  \, (X)$ primitively into $II_{1,25}=
\Lambda\oplus U$ such that the projection of the Weyl vector
$w=(0,(0,1))\in \Lambda\oplus U$ onto $\text{\rm Pic}  \, (X)\otimes{\bold R}$
must be an ample class.

\par \vskip 1pc
Conway [Co] described a fundamental domain $D$ of the reflection
group $W(II_{1,25})^{(2)}$ in terms of the Leech roots (=roots
with intersection number 1 with the Weyl vector $w$). More
precisely, he showed that $W(II_{1,25})^{(2)}$ is generated by
$(-2)$-reflections corresponding to Leech roots.

\par \vskip 1pc \noindent
{\bf Step 3.} The fundamental domain $D$ of the reflection group
$W(II_{1,25})^{(2)}$ cuts out a finite polyhedral cone $D'$ inside
the ample cone $D(X)$. In other words, $$D'=D\cap P_+(X),$$ where
$P_+(X)$ is the positive component of $\{\omega\in
\text{\rm Pic}  \, (X)\otimes{\bold R} : \langle\omega,\omega\rangle>0\}$.
Indeed, $D$ contains the Weyl vector $w$ and, by Step 2, the
projection of $w$ is contained in $D'$.

Determine the hyperplanes $\alpha$ which bound $D'$. The reason
why $D'$ is polyhedral comes from Borcherds [Bor1]; among
infinitely many faces of $D$, those intersecting $P_+(X)$ bound
$D'$, and these faces correspond to Leech roots having a non-zero
projection onto $\text{\rm Pic}  \, (X)\otimes{\bold R}$.

\par \vskip 1pc \noindent
{\bf Step 4.} Match the faces $\alpha$ of $D'$ with automorphisms
$g_\alpha$ such that $g_\alpha$ sends one of the half-spaces
defined by $\alpha$ to the other half-space defined by $\alpha$
or to one of the two half-spaces corresponding to $g_\alpha^{-1}$.

This allows one to prove that the automorphisms $g_\alpha$
generate a group of symmetries of $D(X)$, having $D'$ as its
fundamental domain.

\par \vskip 1pc \noindent
{\bf Step 5.} Take care of $SymD'$. See if which symmetries of
$D'$ realize as automorphisms of $X$. Finally see if there are any
projectively linear automorphisms, (which generate the kernel of
the map in {\bf (5.5)}).

\par \vskip 1pc \noindent
{\bf Remark.} In the known cases {\bf (5.7) - (5.10)}, the
embedding $\text{\rm Pic}  \, (X)\subset \Lambda\oplus U$ is given in such a way
that $\text{\rm Pic}  \, (X)$ is the orthogonal complement of a root sublattice of
$\Lambda\oplus U$. For example, in case {\bf (5.7)} the orthogonal
complement of $\text{\rm Pic}  \, (X)$ in $\Lambda\oplus U$ is a primitive
sublattice of rank 10 which contains a negative definite root
lattice of type $A_5+A_1^5$, which is of index 2.

In practice, Step 4 seems most complicated. The automorphism
$g_\alpha$ works like a reflection, but is not necessarily an
involution. It may be of infinite order. At any rate, $g_\alpha$
may be geometrically evident, or can be picked up from a list of
already known automorphisms, or may be found by looking at extra
structures of $X$, e.g. elliptic fibrations, double plane
structures,..., etc. In worst cases, one has to find an (abstract)
effective Hodge isometry of $H^2(X, {\bold Z})$ and then realize
it geometrically.

The reason why the beautiful combinatorics of the Leech lattice
plays a role in the description of the automorphism groups of K3
surfaces is still unclear to us. We hope that the classification
of all K3 surfaces whose Picard lattice is isomorphic to the
orthogonal complement of a root sublattice of $II_{1,25}$ will
shed more light to this question.

\par \vskip 1pc \noindent
{\bf Remark (5.13).} The method {\bf (5.12)} also works for some
supersingular K3 surfaces. I. Dolgachev and S. Kond$\bar{\roman
o}$ [DoKo] have computed the automorphism group of a supersingular
K3 surface in characteristic 2 whose Picard lattice is $U\oplus
D_{20}$. (For supersingular K3 surfaces Torelli type theorem
holds, i.e. an automorphism of a supersingular K3 surface is
determined by its action on the Picard lattice.)

An even lattice $L$ is reflective if its reflection group $W(L)$
is of finite index in $O(L)$ (see {\bf (5.5)}). The lattice
$U\oplus D_{20}$ is reflective as it is pointed out by Borcherds
[Bor1] and it is the only known example, up to scaling, of an
even reflective hyperbolic lattice of rank 22.  The range of
possible rank of an even reflective lattice of signature
$(1,r-1), r\ge 1$, is given by Esselmann [Es]: it takes the same
range 1, 2,..., 20, 22 as the Picard number of a K3 surface in
positive characteristic.

\head
Section 6. Automorphisms of algebraic surfaces - singular surface case
\endhead

\par \noindent
{\bf (6.1).} In [MM] and [MZ3], one considers pairs $(V, G)$ of surface $V$
and group $G$ of automorphisms, where $V$ may be singular and even non-complete.

\par \vskip 1pc \noindent
{\bf (6.2).} Let $Y$ be a Gorenstein del Pezzo singular surface of Picard number 1.
In [Z9], one classifies all actions on $Y$ by cyclic groups
${\bold Z}/(p)$ of prime order $p \ge 5$.

\par \vskip 1pc \noindent
{\bf Theorem} [Z9, Theorems A and C].

\par \noindent
{\it Let $Y$ be a Gorenstein del Pezzo singular surface
of Picard number $1$. Then we have:}

\par \noindent
(1) {\it Either $|\text{\rm Aut}(Y)| = 2^a3^b$ for some $1 \le a+b \le 7$,
or $\text{\rm Aut}(Y) \supseteq {\bold Z}/(p)$
for every prime $p \ge 5$ and hence $|\text{\rm Aut}(Y)| = \infty$.}

\par \noindent
(2) {\it $|\text{\rm Aut}(Y)|$ is finite if and only if either $\text{\rm Sing} \, Y = A_7$
or $K_Y^2 = 1$ and $|-K_Y|$ has at least three singular members.}

\par \noindent
(3) {\it Let $p \ge 5$ be a prime.
Suppose that $Y$ is not isomorphic to the quadric cone in ${\bold P}^3$.
Then modulo equivariant isomorphism,
there is either none, or only one, or exactly $p+1$ action(s)
of ${\bold Z}/(p)$ on $Y$. All actions are given in $[Z9]$.}

\par \vskip 1pc \noindent
{\bf Remark 6.3.} We like to compare {\bf (6.2)} with known results
for smooth del Pezzo surfaces.

\par \noindent
(1) If $X$ is a generic rational surface with $K_X^2 \le 5$, then
$|\text{\rm Aut}(X)|$ divides $5!$ (see [DO], [Ki]) .

\par \noindent
(2) If $X$ is a del Pezzo surface of degree 3, 4, then $\text{\rm Aut}(X)$ does not
contains ${\bold Z}/(p)$ for any prime $p \ge 7$, and modulo equivariant
isomorphism there is at most one non-trivial ${\bold Z}/(5)$ action on $X$ [Ho1, 2].

\par \noindent
(3) Let $X$ be a rational projective surface with a non-trivial
${\bold Z}/(p)$-action for some prime $p$ such that
the ${\bold Z}/(p)$-invariant sublattice of $\text{\rm Pic}  \,  X$ is of rank 1
(this condition is automatic if the Picard number $\rho(X) = 1$).
If $X$ is smooth, then $p \le 5$ unless $X = {\bold P}^2$
[ZD, Theorem 1].

\head
Acknowledgement\
\endhead
We would like to thank the referee for careful reading and
suggestions which improve the paper.
This article was initiated when both of us were attending the
conference -- Algebraic Geometry in East Asia -- in August 2001,
and we would like to thank the organizers for providing the
wonderful environment and the encouragement for us to write a
survey article. This work was finalized when the second-named
author visited Korea Institute for Advanced Study in December
2001, and he likes to thank the institute for the hospitality.
This work was partially supported by Korea Science and
Engineering Foundation (R01-1999-00004) and an Academic Research
Fund of National University of Singapore.

\head
References
\endhead

\par \noindent
[AN] V. A. Alekseev and V. V. Nikulin,
Classification of del Pezzo surfaces with log-terminal singularities of index
$\le 2$ and involutions on $K3$ surfaces,
Soviet Math. Dokl. {\bf 39} (1989), 525--528.

\par \noindent
[Ar1] M. Artin, Some numerical criteria for contractability of curves
on algebraic surfaces. Amer. J. Math. {\bf 84} (1962), 485--496.

\par \noindent
[Ar2] M. Artin, On isolated rational singularities of surfaces,
Amer. J. Math. {\bf 88} (1966), 129--136.

\par \noindent
[BPV] W. Barth, C. Peters and A. Van de Ven, Compact Complex
Surfaces, Springer-Verlag, 1984.

\par \noindent
[Bl] R. Blache, The structure of l.c. surfaces of Kodaira dimension zero, I,
J. Algebraic Geom. {\bf 4} (1995), 137--179.

\par \noindent
[Bor1] R. Borcherds, Automorphism groups of Lorentzian lattices,
J. Algebra {\bf 111} (1987), 133-153.

\par \noindent
[Bor2] R. Borcherds, Coxter groups, Lorentzian lattices, and K3
surfaces, International Mathematics Research Notices {\bf 19}
(1998), 1011-1031.

\par \noindent
[Br] E. Brieskorn, Rationale Singularit\"{a}ten komplexer
Fl\"{a}chen, Invent. Math. {\bf 4} (1967 / 1968), 336--358.

\par \noindent
[BR] D. Burns, Jr. and M. Rapoport, On the Torelli theorem for
K\"{a}hlerian $K3$ surfaces, Ann. Sci. Ecole Norm. Sup. (4) {\bf
8} (1975), 235-273.

\par \noindent
[Ca] F. Campana, Connexiti rationnelle des variitis de Fano, Ann.
Sci. Ecole Norm. Sup. {\bf 25} (1992), 539--545.

\par \noindent
[CKO] F. Catanese, J. Keum and K. Oguiso,
Some remarks on the universal cover of an open K3 surface,
Math. Ann. to appear, math.AG/{\bf 0104015}.

\par \noindent
[Co] J.H. Conway, The automorphism group of the $26$ dimensional
even Lorentzian lattice, J. Algebra {\bf 80} (1983), 159-163.

\par \noindent
[CD] F. Cossec and I. V. Dolgachev, Enriques surfaces, I,
Progress in Mathematics, {\bf 76} (1989), Birkh\"{a}user Boston,
Inc.

\par \noindent
[dF] Tommaso de Fernex, Birational trasformations of prime order of the projective plane,
preprint 2001.

\par \noindent
[DP] J. -P. Demailly and M. Paun,
Numerical characterization of the K\"ahler cone of a compact K\"ahler manifold,
math.AG / {\bf 0105176}.

\par \noindent
[DK] I. V. Dolgachev and J. Keum, Birational automorphisms of
quartic Hessian surfaces, Trans. Amer. Math. Soc. 354 (2002),
3031-3057.

\par \noindent
[DoKo] I. V. Dolgachev and S. Kond$\bar{\roman o}$, A
supersingular K3 surface in characteristic 2 and the Leech
lattice, math.AG/0112283 v2.

\par \noindent
[DO] I. Dolgachev and D. Ortland, Point sets in projective spaces
and theta functions, Astérisque, Vol. {\bf 165} (1988).

\par \noindent
[DZ] I. V. Dolgachev and D. -Q. Zhang,
Coble rational surfaces, Amer. J. Math. {\bf 123} (2001), 79--114.

\par \noindent
[Du] A. H. Durfee, Fifteen characterizations of rational double points and
simple critical points, Enseign. Math. {\bf 25} (1979), 131--163.

\par \noindent
[Es] F. Esselmann, Ueber die maximale Dimension von
Lorenz-Gittern mit coendlicher Spiegelungsgruppe, J. Number
Theory, {\bf 61} (1996), 103-144.

\par \noindent
[FKL] A. Fujiki, R. Kobayashi and S. Lu,
On the fundamental group of certain open normal surfaces.
Saitama Math. J. {\bf 11} (1993), 15--20.

\par \noindent
[Fu] T. Fujita, On the topology of noncomplete algebraic surfaces,
J. Fac. Sci. Univ. Tokyo Sect. IA Math. {\bf 29} (1982), 503--566.

\par \noindent
[Gi] M. H. Gizatullin, Rational $G$-surfaces,
Math. USSR Izv. {\bf 16} (1981), 103--134.

\par \noindent
[GHS] T. Graber, J. Harris and J. Starr, Families of rationally connected varieties,
math.AG / {\bf 0109220}.

\par \noindent
[GPa] R. V. Gurjar and A.J. Parameswaran,
Affine lines on {\bf Q}-homology planes, J. of Math. Kyoto
Univ. {\bf 35} (1995), 63-77.

\par \noindent
[GPr] R. V. Gurjar and C. R. Pradeep,
$Q$-homology planes are rational, III,
Osaka J. Math. {\bf 36} (1999), 259--335.

\par \noindent
[GPS] R. V. Gurjar, C. R. Pradeep and A. R. Shastri,
On rationality of logarithmic {\bf Q}-homology planes,
II. Osaka J. Math. {\bf 34} (1997), 725--743.

\par \noindent
[GS] R. V. Gurjar and A. R. Shastri,
On the rationality of complex homology $2$-cells, I; II, J. Math.
Soc. Japan, {\bf 41} (1989), 37--56; {\bf 41} (1989), 175--212.

\par \noindent
[GZ1] R. V. Gurjar and D. -Q. Zhang, $\pi_1$ of smooth points of a log del Pezzo
surface is finite, I, J. Math. Sci. Univ. Tokyo {\bf 1} (1994), 137--180.

\par \noindent
[GZ2] R. V. Gurjar and D. -Q. Zhang, $\pi_1$ of smooth points of a log del
Pezzo surface is finite, II, J. Math. Sci. Univ. Tokyo {\bf 2} (1995), 165--196.

\par \noindent
[GZ3] R. V. Gurjar, and D. -Q. Zhang, Normal algebraic surfaces with trivial
bicanonical divisor, J. Algebra {\bf 186} (1996), 970--989.

\par \noindent
[Ho1] T. Hosoh, Automorphism groups of cubic surfaces, J. Algebra {\bf 192} (1997),
651--677.

\par \noindent
[Ho2] T. Hosoh, Automorphism groups of quartic del Pezzo surfaces,
J. Algebra {\bf 185} (1996), 374--389.

\par \noindent
[Hu] J. Hutchinson, The Hessian of the cubic surface, Bul. Amer.
Math. Soc. {\bf 5} (1889), 282-292.

\par \noindent
[I1] S. Iitaka, On logarithmic $K3$ surfaces, Osaka J. Math. {\bf 16} (1979),
675--705.

\par \noindent
[I2] S. Iitaka, A numerical criterion of quasi-abelian surfaces,
Nagoya Math. J. {\bf 73} (1979), 99--115.

\par \noindent
[I3] S. Iitaka, Algebraic geometry -- An introduction to birational geometry
of algebraic varieties, Graduate Texts in Mathematics, {\bf 76},
Springer-Verlag, New York-Berlin, 1982.

\par \noindent
[Is] V. A. Iskovskih, Minimal models of rational surfaces over
arbitrary fields, Math. USSR Izv. {\bf 14} (1981), 17--39.

\par \noindent
[Kt] S. Kantor, Theorie der endlichen Gruppen von eindeutigen Transformationen
in der Ebene, Berlin: Mayer $\&$ Muller, 1895.

\par \noindent
[KN] M. Kato and I. Naruki,
Depth of rational double points on quartic surfaces,
Proc. Japan Acad. Ser. A Math. Sci. 58 (1982), 72--75.

\par \noindent
[Ka1] Y. Kawamata, On the classification of noncomplete algebraic
surfaces, Lecture Notes in Math. {\bf 732} (1979), 215--232.
Springer, Berlin.

\par \noindent
[Ka2] Y. Kawamata, The cone of curves of algebraic varieties,
Ann. of Math. {\bf 119} (1984), 603--633.

\par \noindent
[KMM] Y. Kawamata, K. Matsuda and K. Matsuki,
Introduction to the minimal model problem,
Adv. Stud. Pure Math. {\bf 10} (1987), 283--360.

\par \noindent
[KMc] S. Keel and J. McKernan, Rational curves on quasi-projective surfaces,
Mem. Amer. Math. Soc. {\bf 140} (1999), no. 669.

\par \noindent
[KMo] J. Kollar and S. Mori, Birational geometry of algebraic varieties,
Cambridge Tracts in Mathematics, {\bf 134} (1998).

\par \noindent
[Ke1] J. Keum, Automorphisms of Jacobian Kummer surfaces,
Compositio Math. {\bf 107} (1997), 269-288.

\par \noindent
[Ke2] J. Keum,  Every algebraic Kummer surface has infinitely
many automorphisms, unpublished manuscript (1996).

\par \noindent
[Ke3] J. Keum, Automorphisms of a generic Jacobian Kummer
surface, Geom. Ded. {\bf 76} (1999), 177-181.

\par \noindent
[KK] J. Keum and S. Kondo,
The automorphism groups of Kummer surfaces associated
with the product of two elliptic curves, Trans. Amer. Math. Soc. {\bf 353} (2001),
1469--1487.

\par \noindent
[KZ] J. Keum and D. -Q. Zhang,
Fundamental groups of open K3 surfaces, Enriques surfaces and Fano 3-folds,
J. Pure App. Alg. {\bf 170} (2002), 67--91.

\par \noindent
[Kb1] R. Kobayashi, Uniformization of complex surfaces. K\"{a}hler
metric and moduli spaces, Adv. Stud. Pure Math. {\bf 18-II}
(1990), 313--394.

\par \noindent
[Kb2] R. Kobayashi, Einstein-K\"{a}hler $V$-metrics on open Satake
$V$-surfaces with isolated quotient singularities, Math. Ann.
{\bf 272} (1985), 385--398.

\par \noindent
[KNS] R. Kobayashi, S. Nakamura and F. Sakai,
A numerical characterization of ball
quotients for normal surfaces with branch loci,
Proc. Japan Acad. Ser. A Math. Sci. {\bf 65} (1989), 238--241.

\par \noindent
[Ki] M. Koitabashi, Automorphism groups of generic rational surfaces,
J. Algebra {\bf 116} (1988), 130 --142.

\par \noindent
[Kj1] H. Kojima, Open rational surfaces with logarithmic Kodaira dimension zero,
Internat. J. Math. {\bf 10} (1999), 619--642.

\par \noindent
[Kj2] H. Kojima, Logarithmic del Pezzo surfaces of
rank one with unique singular points,
Japan. J. Math. (N.S.) {\bf 25} (1999), 343--375.

\par \noindent
[Kj3] H. Kojima, Complements of plane curves with logarithmic Kodaira dimension zero,
J. Math. Soc. Japan {\bf 52} (2000), 793--806.

\par \noindent
[Kj4] H. Kojima, Rank one log del Pezzo surfaces of index two, J. Math. Kyoto Univ.
to appear.

\par \noindent
[KoMiMo] J. Kollar, Y. Miyaoka and S. Mori,
Rationally connected varieties, J. Alg. Geom. {\bf 1} (1992), 429--448.

\par \noindent
[Kon1] S. Kond$\bar{\roman o}$, The automorphism group of a
generic Jacobian Kummer surface, J. Alg. Geom. {\bf 7}
(1998), 589-609.

\par \noindent
[Kon2] S. Kond$\bar{\roman o}$, Niemeier lattices, Mathieu
groups, and finite groups of symplectic automorphisms of K3
surfaces (with an appendix by S. Mukai), Duke Math. J. {\bf 92}
(1998), 593-603.

\par \noindent
[Kon3] S. Kond$\bar{\roman o}$, The maximum order of finite
groups of automorphisms of $K3$ surfaces, Amer. J. Math.
{\bf 121} (1999), 1245-1252.

\par \noindent
[Kon4] S. Kond$\bar{\roman o}$, Enriques surfaces with finite
automorphism groups, Japan. J. Math. {\bf 12} (1986),
191-282.

\par \noindent
[Ma1] Yu. I. Manin, Rational surfaces over perfect fields, II, Math. USSR Sb.
{\bf 1} (1967), 141--168.

\par \noindent
[Ma2] Yu. I. Manin, Cubic forms : Algebra, geometry, arithmetic,
2nd ed, North-Holland Math. Library, {\bf 4} (1986),
North-Holland Publ. Co., Amsterdam-New York.

\par \noindent
[Mi] M. Miyanishi, Open algebraic surfaces, CRM Monograph Series, {\bf 12},
American Mathematical Society, 2001.

\par \noindent
[MM] M. Miyanishi and K. Masuda, Open algebraic surfaces with finite group actions,
Transform. Group, to appear.

\par \noindent
[MT1] M. Miyanishi and S. Tsunoda,
Noncomplete algebraic surfaces with logarithmic
Kodaira dimension $-\infty$ and with nonconnected boundaries at infinity,
Japan. J. Math. (N.S.) {\bf 10} (1984), 195--242.

\par \noindent
[MT2] M. Miyanishi and S. Tsunoda,
Logarithmic del Pezzo surfaces of rank one with
noncontractible boundaries, Japan. J. Math. (N.S.) {\bf 10} (1984), 271--319.

\par \noindent
[MZ1] M. Miyanishi and D. -Q. Zhang, Gorenstein log del Pezzo surfaces of rank one,
J. Algebra {\bf 118} (1988), 63--84.

\par \noindent
[MZ2] M. Miyanishi and D. -Q. Zhang,
Gorenstein log del Pezzo surfaces, II, J. Algebra {\bf 156} (1993), 183--193.

\par \noindent
[MZ3] M. Miyanishi and D. -Q. Zhang, Equivariant classification of Gorenstein open
log del Pezzo surfaces with finite group actions, Preprint 2001.

\par \noindent
[Mu] S. Mukai, Finite groups of automorphisms of $K3$ surfaces
and the Mathieu group, Invent. Math. {\bf 94} (1988), 183-221.

\par \noindent
[Ni1] V.V. Nikulin, On the quotient groups of the automorphism
group of hyperbolic forms by the subgroups generated by
2-reflections, J. Soviet Math. {\bf 22} (1983), 1401-1476.

\par \noindent
[Ni2] V.V. Nikulin, On the description of groups of automorphisms
of Enriques surfaces, Soviet Math. Dokl. {\bf 277} (1984),
1324-1327.

\par \noindent
[Ni3] V.V. Nikulin, Del Pezzo surfaces with log-terminal
singularities, I; II; III, Math. USSR-Sb. {\bf 66} (1990),
231--248; Math. USSR-Izv. {\bf 33} (1989), 355--372; Math.
USSR-Izv. {\bf 35} (1990), 657--675.

\par \noindent
[Ni4] V. V. Nikulin, Discrete reflection groups in Lobachevsky spaces and algebraic surfaces, Proc. Intern. Cong. Math. Berkeley, Calif. 1986, pp.654--671.

\par \noindent
[Og1] K. Oguiso, Picard numbers in a family of hyperk\"ahler manifolds
and applications, J. Alg. Geom. to appear

\par \noindent
[Og2] K. Oguiso, Picard numbers in a family of hyperk\"ahler manifolds -
A supplement to the article of R. Borcherds, L. Katzarkov, T. Pantev,
N. I. Shepherd-Barron, math.AG / {\bf 0011258}.

\par \noindent
[Og3] K. Oguiso, Automorphism groups in a family of K3 surfaces,
math.AG / {\bf 0104049}.

\par \noindent
[OZ1] K. Oguiso and D. -Q. Zhang, On the most algebraic $K3$ surfaces and
the most extremal log Enriques surfaces, Amer. J. Math. {\bf 118} (1996),
1277--1297.

\par \noindent
[OZ2] K. Oguiso and D. -Q. Zhang, On the complete classification of
extremal log Enriques surfaces, Math. Z. {\bf 231} (1999), 23--50.

\par \noindent
[OZ3] K. Oguiso and D. -Q. Zhang, K3 surfaces with order 11 automorphisms,
math.AG / {\bf 9907020}.

\par \noindent
[OZ4] K. Oguiso and D. -Q. Zhang, On Vorontsov's theorem on $K3$ surfaces
with non-symplectic group actions, Proc. Amer. Math. Soc. {\bf 128} (2000),
1571--1580.

\par \noindent
[OZ5] K. Oguiso and D. -Q. Zhang, The Simple Group of Order 168 and K3 Surfaces,
in : Complex Geometry, Collection of papers dedicated to Hans Grauert;
Bauer, Catanese, Kawamata, Peternell and Siu (ed), Springer, 2002;
math.AG / {\bf 0011259}.

\par \noindent
[Oh] K. Ohno, Toward determination of the singular fibers of
minimal degeneration of surfaces with $\kappa = 0$, Osaka J. Math. {\bf 33}
(1996), 235--305.

\par \noindent
[PSS] I. Piatetski-Shapiro, I.R. Shafarevich, A Torelli theorem
for algebraic surfaces of type K3, Math. USSR Izv. {\bf 5}
(1971), 547-587.

\par \noindent
[PS] C. R. Pradeep and A. R. Shastri,
On rationality of logarithmic {\bf Q}-homology planes, I,
Osaka J. Math. {\bf 34} (1997), 429--456.

\par \noindent
[Re1] M. Reid, Young person's guide to canonical singularities,
Proc. Sympos. Pure Math. Part 1 {\bf 46} (1987), 345--414.

\par \noindent
[Re2] M. Reid, Campedelli versus Godeaux --
Problems in the theory of surfaces and their
classification, Sympos. Math., XXXII (1991), 309--365.

\par \noindent
[Sa] F. Sakai, Semistable curves on algebraic surfaces and
logarithmic pluricanonical maps, Math. Ann. {\bf 254} (1980), 89--120.

\par \noindent
[Se] B. Segre, The Non-singular cubic surfaces,
Oxford University Press, Oxford, 1942.

\par \noindent
[SZ] I. Shimada and D. -Q. Zhang, Classification of extremal elliptic
$K3$ surfaces and fundamental groups of open $K3$ surfaces, Nagoya Math. J.
{\bf 161} (2001), 23--54.

\par \noindent
[Su] K. Suzuki, On Morrison's cone conjecture for Klt surfaces with
$K_X \equiv 0$, Comment. Math. Univ. St. Paul {\bf 50} (2001), 173--180.

\par \noindent
[Siu] Y. T. Siu,  Every $K3$ surface is K\"{a}hler, Invent. Math.
{\bf 73} (1983), 139-150.

\par \noindent
[Ta] S. Takayama, Simple connectedness of weak Fano varieties,
J. Algebraic Geom. {\bf 9} (2000), 403--407.

\par \noindent
[Ts] S. Tsunoda, Structure of open algebraic surfaces,
I, J. Math. Kyoto Univ. {\bf 23} (1983), 95--125.

\par \noindent
[TZ] S. Tsunoda and D. -Q. Zhang,
Noether's inequality for noncomplete algebraic surfaces of
general type, Publ. Res. Inst. Math. Sci. {\bf 28} (1992), 21--38.

\par \noindent
[Vin] E.B. Vinberg, The two most algebraic K3 surfaces, Math.
Ann. {\bf 265} (1985), 1-21.

\par \noindent
[Xi1] G. Xiao, Bound of automorphisms of surfaces of general type.
I, Ann. of Math. (2) {\bf 139} (1994), no. 1, 51-77.

\par \noindent
[Xi2] G. Xiao, Bound of automorphisms of surfaces of general type.
II, J. Algebraic Geom. {\bf 4} (1995), no. 4, 701-793.

\par \noindent
[Z1] D. -Q. Zhang, On Iitaka surfaces, Osaka J. Math. {\bf 24} (1987), 417--460.

\par \noindent
[Z2] D. -Q. Zhang, Logarithmic del Pezzo surfaces of rank one
with contractible boundaries, Osaka J. Math. {\bf 25} (1988), 461--497.

\par \noindent
[Z3] D. -Q. Zhang, Logarithmic del Pezzo surfaces with rational double and
triple singular points, Tohoku Math. J. {\bf 41} (1989), 399--452.

\par \noindent
[Z4] D. -Q. Zhang, Logarithmic Enriques surfaces, I; II,
J. Math. Kyoto Univ. {\bf 31} (1991), 419--466; {\bf 33} (1993), 357-397.

\par \noindent
[Z5] D. -Q. Zhang, Noether's inequality for noncomplete algebraic surfaces
of general type, II, Publ. Res. Inst. Math. Sci. {\bf 28} (1992), 679--707.

\par \noindent
[Z6] D. -Q. Zhang, Algebraic surfaces with nef and big anti-canonical divisor,
Math. Proc. Cambridge Philos. Soc. {\bf 117} (1995), 161--163.

\par \noindent
[Z7] D. -Q. Zhang, The fundamental group of the smooth part of a log Fano variety,
Osaka J. Math. {\bf 32} (1995), 637--644.

\par \noindent
[Z8] D. -Q. Zhang, Algebraic surfaces with log canonical singularities and
the fundamental groups of their smooth parts, Trans. Amer. Math. Soc.
{\bf 348} (1996), 4175--4184.

\par \noindent
[Z9] D. -Q. Zhang, Automorphisms of finite order on Gorenstein del Pezzo surfaces,
Trans. Amer. Math. Soc. {\bf 354} (2002), 4831--4845.

\par \noindent
[ZD] D. -Q. Zhang, Automorphisms of finite order on rational surfaces,
with an appendix by I. Dolgachev, J. Algebra {\bf 238} (2001),560--589.

\par \vskip 3pc \noindent
J. Keum
\par \noindent
Korea Institute for Advanced Study
\par \noindent
207-43 Cheongryangri-dong, Dongdaemun-gu
\par \noindent
Seoul 130-012, Korea
\par \noindent
E-mail : jhkeum$\@$kias.re.kr

\par \vskip 1pc \noindent
D. -Q. Zhang
\par \noindent
Department of Mathematics
\par \noindent
National University of Singapore
\par \noindent
2 Science Drive 2, Singapore 117543
\par \noindent
Republic of Singapore
\par \noindent
E-mail: matzdq$\@$math.nus.edu.sg
\enddocument